\documentclass[12pt]{article}
\usepackage[utf8]{inputenc}
\usepackage[T1]{fontenc}
\usepackage{amsmath, amssymb, amsthm}
\usepackage{geometry}
\usepackage{graphicx}
\usepackage{enumitem}
\usepackage{authblk}
\usepackage{url}
\usepackage{xcolor}

\newtheorem{theorem}{Theorem}

\newtheorem{conjecture}[theorem]{Conjecture}

\title{{Symplectic Widths and Filters in Classical Reaction Dynamics:\\
       Normal-Form Bottlenecks, Energy Layers, and Finite-Time Diagnostics}}
\author[1]{Stephen Wiggins\thanks{stephen.wiggins@me.com}}
\affil[1]{Hetao Institute of Mathematics and Interdisciplinary Sciences, Shenzhen, China}
\affil[2]{School of Mathematics, University of Bristol, UK}
\date{\today}

\begin{document}

\maketitle

\begin{abstract}
We develop a symplectic-geometric framework for studying classical reaction dynamics near an index-1 saddle.
The central idea is that a reaction bottleneck can be viewed as a \emph{symplectic filter}: a transition-state region whose transverse compressibility is constrained not only by energy and flux, but also by Gromov's non-squeezing theorem and the associated notion of symplectic capacity.
Using Poincar\'e--Birkhoff normal form theory, we describe the phase-space structures organizing transport through the saddle region, including dividing surfaces, normally hyperbolic invariant manifolds (NHIMs), and the bath-action geometry of the transverse modes.

A recurring issue is that a fixed-energy surface is odd-dimensional and is not itself the right object to which a symplectic capacity can be assigned.
We therefore distinguish carefully between three related objects: the fixed-energy
transition-state geometry, a full-dimensional energy layer obtained by thickening
the energy surface, and explicit phase-space ensembles such as balls or
bath-localized ensembles.
For quadratic saddle--center and saddle--center--center models, the normal-form geometry identifies natural transverse bath-action area scales associated with the bottleneck.
For anharmonic Eckart--Morse and Eckart--Morse--Morse normal forms, we formulate corresponding \emph{candidate} symplectic width scales based on maximal admissible bath actions in a bounded local neighborhood of the saddle.

The numerical experiments are presented as diagnostics, not as proofs of a universal reaction-rate theorem.
They illustrate how finite ensembles can exhibit transverse spreading, bath-action localization, and finite-time delay in ways that are not resolved by scalar diagnostics such as total Liouville volume or total flux alone.
The result is a capacity-based geometric framework for studying reaction
bottlenecks, linking normal-form transition-state geometry with
finite-time diagnostics of bath-action localization and transverse
phase-space spreading.
\end{abstract}


\section{Introduction}

The phase-space formulation of transition state theory (TST) has revealed a rich geometric structure underlying chemical reactions: dividing surfaces with minimal-flux and no-local-recrossing properties, normally hyperbolic invariant manifolds (NHIMs), stable and unstable manifolds that organize transport, and the associated gap-time and flux descriptions of reaction dynamics.
These structures are now computable for multidimensional molecular models using Poincar\'e--Birkhoff normal form theory, and they have been used to compute reaction rates, reactive volumes, gap-time distributions, and other dynamical quantities with high accuracy \cite{Waalkens2008,Uzer2002,Waalkens2005a,Waalkens2005b}.

The present paper adds a further geometric ingredient to this phase-space picture.
Gromov's non-squeezing theorem \cite{Gromov1985} shows that Hamiltonian flows are more restrictive than general volume-preserving flows.
Liouville's theorem says that Hamiltonian flow preserves full phase-space volume.
Non-squeezing says something stronger and more directional: a symplectic image of a ball cannot be embedded into a canonical cylinder whose coordinate--momentum cross-section is too small.
The associated invariant is the symplectic capacity, introduced in the work of Ekeland and Hofer \cite{EkelandHofer1989,EkelandHofer1990}.
For a ball of radius $R$, this capacity is $\pi R^2$; for a canonical cylinder of radius $r$, it is $\pi r^2$.
Thus the relevant comparison is not ordinary $2n$-dimensional volume, but a two-dimensional symplectic size measured in canonical conjugate planes.

Our guiding interpretation is that a reaction bottleneck can act as a \emph{symplectic filter}.
The word ``filter'' is meant to emphasize that passage through the transition-state region is controlled not only by total energy or by the total flux through a dividing surface, but also by how a finite phase-space set is arranged relative to the canonical transverse bath directions.
In normal-form coordinates near an index-1 saddle, one pair of variables describes the reaction direction and the remaining pairs describe bounded transverse, or bath, modes.
The NHIM bounds the dividing surface and determines natural bath-action scales.
These scales suggest canonical cylinders against which one can compare the symplectic size of incoming phase-space sets.

It is essential, however, not to overstate what Gromov's theorem by itself proves.
The theorem is a statement about embeddings of full-dimensional symplectic domains.
It does not identify which individual trajectory reacts, returns, or remains delayed.
It also does not compute a reaction probability or a reaction rate by itself.
Rather, it supplies a topological obstruction to coherent squeezing of a full phase-space set into a smaller canonical cylinder.
The dynamical consequences of that obstruction---transverse spreading, bath-action redistribution, incomplete finite-time transmission, or global return---must be studied using the actual Hamiltonian flow.
This distinction between theorem-level non-containment and model-dependent trajectory dynamics is central to the interpretation developed below.

A second point is equally important for reaction dynamics.
Let $z$ denote a point in the $2n$-dimensional phase space, and let
$H(z)$ be the Hamiltonian function.
The traditional microcanonical setting is the fixed-energy surface
\[
        \Sigma_E=\{z:H(z)=E\},
\]
where $E$ is the chosen total energy.
This is the natural setting for many formulations of classical reaction dynamics.
However, $\Sigma_E$ is a $(2n-1)$-dimensional hypersurface inside the full
phase space, and hence it is not itself the kind of full-dimensional
symplectic domain to which Gromov's non-squeezing theorem applies directly.
Therefore one should not assign a finite symplectic capacity directly to
the fixed-energy surface alone.

To pose a non-squeezing question near a fixed energy, one must instead work
with a full-dimensional domain in the ambient phase space.
There are two natural ways to do this.
One may use an explicit full-dimensional ensemble, such as a ball, ellipsoid,
or other localized phase-space set.
Alternatively, one may thicken the fixed-energy surface to a narrow energy layer
\[
        \mathcal L_{E,\Delta E}
        =
        \{z:E-\Delta E\le H(z)\le E+\Delta E\},
\]
where $\Delta E>0$ is a small energy thickness.
This layer has full dimension in phase space.
It may then be further localized near the transition-state region by using
flow-time or gap-time coordinates.
In this interpretation, $E$ is the central energy of the layer and is the
energy at which the leading-order bottleneck geometry and candidate width
scales are evaluated.

The paper therefore has two connected aims.
The first is to identify natural transverse symplectic width scales associated with the normal-form bottleneck.
For quadratic saddle--center--$\cdots$--center models these scales can be written explicitly in terms of the maximal bath actions on the NHIM.
For anharmonic normal forms, we define analogous candidate scales by solving for the maximal admissible bath actions in the truncated normal-form Hamiltonian.
The second aim is to examine how these scales may appear in finite-time trajectory diagnostics.
The numerical experiments below are exploratory: they are intended to show how transverse localization in bath-action directions can produce strong finite-time effects, not to prove a universal theorem about reaction rates.

The paper is organized as follows.
Section 2 recalls the symplectic structure of Hamiltonian mechanics and the distinction between symplectic maps, volume preservation, and projection areas.
Section 3 reviews Gromov's non-squeezing theorem, symplectic capacity, and the capacity of ellipsoids.
Section 4 recalls the normal-form phase-space structures near saddle--center--$\cdots$--center equilibria.
Section 5 formulates candidate symplectic width scales for quadratic and anharmonic bottlenecks, emphasizing the need for full-dimensional proxy domains or energy layers.
Section 6 describes the Eckart--Morse normal-form models used in the computations.
Section 7 presents the numerical diagnostics.
The concluding section summarizes the symplectic-filter interpretation and the open mathematical questions.

\subsection{Scope and limitations}

The scope of the paper is deliberately classical.
Although symplectic capacity is sometimes discussed in relation to uncertainty principles and quantum mechanics, no quantum dynamics is used here.
The relevant objects are classical Hamiltonian flows, classical normal forms, and classical phase-space ensembles.

The main limitation is also explicit.
We do not claim to have proved that the action-based quantities introduced below are the exact Gromov widths of a uniquely defined reactive domain for an arbitrary molecular Hamiltonian.
Instead, we use normal-form geometry to construct physically meaningful candidate width scales for bounded, full-dimensional neighborhoods of the bottleneck.
In the quadratic model these quantities have a direct canonical-cylinder interpretation.  In the anharmonic normal-form setting they are geometrically well motivated, but the identification of the exact symplectic capacity of a fully specified reactive neighborhood remains an open problem.

The systems studied are:
\begin{itemize}
    \item Quadratic saddle--center and saddle--center--center models, for which the normal-form geometry is explicit and the bath-action width scales can be computed analytically.
    \item Eckart--Morse and Eckart--Morse--Morse models, for which high-order classical normal forms have been computed previously \cite{Waalkens2008,Waalkens2004}.  These models provide a controlled anharmonic setting in which to test the diagnostic value of the candidate width scales.
\end{itemize}

\section{Hamiltonian Mechanics and Symplectic Geometry}

In this section, we briefly establish the mathematical terminology and notation that will be used throughout the paper.
We first introduce the symplectic form, which provides the fundamental geometric structure of phase space, and then define symplectic transformations.
This geometric machinery is necessary to properly state Gromov's non-squeezing theorem in Section 3 and to understand the projection-area constraints evaluated in our numerical experiments in Section 7.

\subsection{Phase space and the symplectic form}

For a comprehensive treatment of the geometric formulation of Hamiltonian mechanics, we refer the reader to the foundational texts by Arnold \cite{Arnold1989} and Abraham and Marsden \cite{AbrahamMarsden1978}.
For a system with \(n\) degrees of freedom, the phase space is \(\mathbb{R}^{2n}\) with coordinates \(z = (q_1,\dots,q_n,p_1,\dots,p_n)\).
The canonical symplectic form is
\[
\omega = \sum_{i=1}^n dp_i \wedge dq_i .
\]
Here, the wedge product (\(\wedge\)) ensures the antisymmetric property of the form, which algebraically encodes the concept of oriented area in the conjugate \((q_i, p_i)\) planes.
A Hamiltonian function \(H(z)\) generates the Hamiltonian vector field \(X_H\) through the relation \(\iota_{X_H}\omega = dH\), where \(\iota\) denotes the interior product.
While this coordinate-free differential geometry notation is elegant, its physical meaning is simply the definition of Hamilton's equations.
In the standard coordinate representation, this relation immediately yields the familiar matrix form \(\dot{z} = J\nabla H(z)\), where the standard symplectic matrix is
\[
J = \begin{pmatrix} 0 & I_n \\ -I_n & 0 \end{pmatrix}.
\]

\subsection{Symplectic maps and integral invariants}

A diffeomorphism \(f\) is \emph{symplectic} (or canonical) if \(f^*\omega = \omega\).
Here, the superscript \(*\) denotes the pullback, which is the formal geometric way of stating that the transformation preserves the symplectic structure.
In coordinates, this is equivalent to the condition that the Jacobian matrix \(Df(z)\) of the transformation satisfies \(Df(z)^{\mathsf T} J Df(z) = J\) everywhere.
This symplectic condition implies several fundamental conservation laws:
\begin{itemize}
    \item Volume preservation (Liouville's theorem): \(\det(Df)=1\), where \(\det\) denotes the determinant of the Jacobian matrix.
    \item Preservation of Poincar\'e integral invariants: for any arbitrary closed loop \(\gamma\) anywhere in the full \(2n\)-dimensional phase space (it need not be restricted to a specific energy surface or conjugate plane), the action integral is conserved:
    \[
    \oint_\gamma \sum p_i dq_i = \oint_{f(\gamma)} \sum p_i dq_i .
    \]
\end{itemize}

A \emph{conjugate plane} (or \emph{symplectic plane}) is a 2-dimensional subspace spanned by a coordinate \(q_i\) and its conjugate momentum \(p_i\);
the restriction of \(\omega\) to such a plane is the standard area form.
It is a crucial mathematical nuance that for linear symplectic maps, preservation of the full \(2n\)-dimensional symplectic form does not in general imply the preservation of ordinary orthogonal projection area onto a fixed conjugate plane.
For nonlinear maps, the area of a projection is not preserved pointwise, though the integral of the symplectic form over any oriented 2-dimensional surface is invariant \cite{Arnold1989, deGosson2006}.
This distinction between preservation of \(2n\)-dimensional phase-space volume and the non-invariance of ordinary two-dimensional projection areas is precisely the setting in which Gromov's non-squeezing theorem gives additional information.

\section{Gromov's Non-Squeezing Theorem and Symplectic Capacity}

We now recall the symplectic-topological result that motivates the rest of the paper.
The purpose of this section is not to give a full introduction to symplectic topology, but to isolate the pieces needed for reaction dynamics: non-squeezing, symplectic capacity, the distinction between capacity and projected area, and the computation of capacities for ellipsoids.

\subsection{Motivation and the symplectic camel}

Liouville's theorem says that Hamiltonian flow preserves $2n$-dimensional phase-space volume.
Volume preservation is fundamental, but it is not the whole symplectic story.
A general volume-preserving map can stretch a set into a very long and very thin filament, making its cross-section in a chosen two-dimensional direction arbitrarily small while keeping the total volume fixed.
Hamiltonian flows are more restrictive: they preserve the symplectic form, and this imposes constraints that cannot be seen from volume alone.

Gromov's non-squeezing theorem \cite{Gromov1985} is the clearest expression of this extra rigidity.
Arnold popularized the result through the ``symplectic camel'' metaphor \cite{Arnold1986}: the phase-space ball is the camel and the narrow canonical cylinder is the eye of the needle.
The crucial word is canonical.
The obstruction concerns containment in cylinders whose bases lie in
canonical coordinate--momentum planes; it is not a statement about arbitrary
two-dimensional projections.
This is the feature that makes the theorem relevant to reaction dynamics, because normal-form coordinates near an index-1 saddle identify natural canonical bath planes transverse to the reaction direction.

\subsection{The non-squeezing theorem}

We now state the form of Gromov's non-squeezing theorem used in this
paper.  Let
\[
        z=(q_1,\ldots,q_n,p_1,\ldots,p_n)\in \mathbb R^{2n}
\]
denote canonical coordinates on phase space.  The canonical symplectic form
is
\[
        \omega=\sum_{j=1}^n dp_j\wedge dq_j .
\]
Equivalently, in matrix notation, the symplectic structure is represented
by the standard matrix
\[
        J=
        \begin{pmatrix}
        0 & I_n\\
        -I_n & 0
        \end{pmatrix}.
\]

For $r>0$, define the $2n$-dimensional ball of radius $r$ by
\[
        B^{2n}(r)
        =
        \left\{
        z\in\mathbb R^{2n}:
        \sum_{j=1}^n(q_j^2+p_j^2)\le r^2
        \right\}.
\]
For $R>0$, define the canonical cylinder of radius $R$, based on the
coordinate--momentum plane $(q_1,p_1)$, by
\[
        Z^{2n}(R)
        =
        \left\{
        z\in\mathbb R^{2n}: q_1^2+p_1^2\le R^2
        \right\}.
\]
The remaining coordinates
\[
        q_2,\ldots,q_n,p_2,\ldots,p_n
\]
are unrestricted.  Thus $Z^{2n}(R)$ is the product of a disk of radius $R$
in the canonical plane $(q_1,p_1)$ with the full space in all other
coordinates.

A symplectic embedding of $B^{2n}(r)$ into $Z^{2n}(R)$ is a smooth injective
map
\[
        f:B^{2n}(r)\to Z^{2n}(R)
\]
that preserves the symplectic structure.  In coordinates, this means that
the derivative matrix $Df(z)$ satisfies
\[
        Df(z)^T J Df(z)=J
\]
at every point $z$ of the ball.  In differential-geometric notation the same
condition is written
\[
        f^\ast\omega=\omega,
\]
where $f^\ast\omega$ denotes the pullback of the two-form $\omega$ by the
map $f$.  Thus a symplectic embedding preserves the Hamiltonian phase-space
structure, not merely the $2n$-dimensional volume.

Gromov's non-squeezing theorem states that if there exists a symplectic
embedding
\[
        f:B^{2n}(r)\to Z^{2n}(R),
\]
then necessarily
\[
        r\le R .
\]
Equivalently, a ball of radius $r$ cannot be embedded symplectically into a
canonical cylinder of smaller radius $R<r$ \cite{Gromov1985}.  This is the
non-squeezing obstruction.

The theorem is substantially stronger than Liouville's theorem.  Liouville's
theorem says that Hamiltonian flow preserves $2n$-dimensional phase-space
volume.  Volume preservation alone would not prevent a set from being
stretched into a long, thin shape and placed inside a cylinder with a very
small two-dimensional base.  Gromov's theorem rules this out when the map is
symplectic and the cylinder is based on a canonical coordinate--momentum
plane.

The restriction to a canonical coordinate--momentum plane is essential.  The
obstruction concerns containment in cylinders whose bases lie in symplectic
planes such as $(q_j,p_j)$.  It is not a statement about arbitrary
two-dimensional projections, nor about cylinders based on non-canonical
coordinate planes such as $(q_1,q_2)$.

For Hamiltonian dynamics, the relevance is immediate.  The time-$t$ map
$\Phi^t$ generated by Hamilton's equations is symplectic wherever it is
defined.  Therefore, if an initial phase-space ball $B^{2n}(r)$ is evolved
by Hamiltonian flow, its image $\Phi^t(B^{2n}(r))$ cannot be contained in a
canonical cylinder $Z^{2n}(R)$ with $R<r$.  This is the precise topological
constraint that will be compared below with the transverse coordinate--
momentum cylinders associated with reaction bottlenecks.

\subsection{Symplectic capacity and projected area}

The theorem motivates the notion of symplectic capacity.
A symplectic capacity is a function $c$ assigning to suitable subsets $\Omega\subset\mathbb R^{2n}$ a number $c(\Omega)\in[0,\infty]$ satisfying the standard axioms of monotonicity, conformality, symplectic invariance, and normalization \cite{EkelandHofer1989,EkelandHofer1990,HoferZehnder1994,McDuffSalamon1998,deGosson2006}:
\begin{enumerate}[label=(\roman*)]
    \item if $\Omega_1\subset\Omega_2$, then $c(\Omega_1)\le c(\Omega_2)$;
    \item $c(\lambda\Omega)=\lambda^2c(\Omega)$ for $\lambda>0$;
    \item $c(f(\Omega))=c(\Omega)$ for every symplectomorphism $f$;
    \item $c(B^{2n}(r))=\pi r^2$ and $c(Z^{2n}(r))=\pi r^2$.
\end{enumerate}
The last normalization is the bridge to non-squeezing: the capacity of a ball and the capacity of a canonical cylinder are both measured by the area of the corresponding two-dimensional disk.

For the present application it is important to distinguish three notions of size.
First, there is ordinary $2n$-dimensional volume, the quantity preserved by Liouville's theorem.
Second, there is the ordinary area of a projection onto a chosen two-dimensional plane.
For an evolving set $\Phi^t(\Omega)$, its projected area in a bath plane might be written as
\[
        A_{q_kp_k}(t)
        =
        \operatorname{Area}\left(\operatorname{proj}_{(q_k,p_k)}\Phi^t(\Omega)\right).
\]
This projected area is a visible diagnostic of deformation, but it is not a symplectic invariant.
It can grow, shrink, or oscillate under Hamiltonian evolution.
Third, there is symplectic capacity, which is invariant under Hamiltonian flow and monotone under symplectic embeddings.

The source of possible confusion is that for an initial round ball $B^{2n}(R)$, the projection onto any canonical coordinate--momentum plane is a disk of area $\pi R^2$, and this is also the ball's capacity:
\[
        A_{q_kp_k}(0)=\pi R^2=c(B^{2n}(R)).
\]
This equality is special to the initial round ball.
After Hamiltonian evolution,
\[
        c(\Phi^t(B^{2n}(R)))=\pi R^2,
\]
but the projected area $A_{q_kp_k}(t)$ need not remain equal to $\pi R^2$.
Thus when numerical experiments below plot projection areas or convex-hull areas, those curves should be read as diagnostics of transverse spreading, not as plots of capacity itself.
The capacity is the invariant obstruction; projected area is a non-invariant
two-dimensional diagnostic used to monitor how the set deforms while
respecting that obstruction.

Two extremal capacities are often useful.
The Gromov width is
\[
        c_{\min}(\Omega)
        =
        \sup\{\pi r^2: B^{2n}(r) \text{ symplectically embeds into } \Omega\},
\]
and the cylindrical capacity is
\[
        c_{\max}(\Omega)
        =
        \inf\{\pi R^2: \Omega \text{ symplectically embeds into } Z^{2n}(R)\}.
\]
For any symplectic capacity one has the usual comparison
\[
        c_{\min}(\Omega)\le c(\Omega)\le c_{\max}(\Omega).
\]
For general domains these quantities can be difficult to compute.
This is why the later formulas in this paper are described as candidate width scales unless a specific domain is simple enough for its capacity to be known exactly.

\subsection{Symplectic spectrum and the capacity of an ellipsoid}

There is one important class of bounded domains for which capacities can be computed explicitly: ellipsoids.
Let $M$ be a real symmetric positive-definite $2n\times 2n$ matrix and consider
\[
        \Omega_M=\{z\in\mathbb R^{2n}:z^{\mathsf T}Mz\le 1\}.
\]
The ordinary eigenvalues of $M$ describe the ellipsoid under orthogonal transformations.
Hamiltonian mechanics instead singles out symplectic transformations, and the relevant invariants are the symplectic eigenvalues.
Let $J$ denote the standard symplectic matrix.
The eigenvalues of $JM$ occur in pairs $\pm i\lambda_j$, with $\lambda_j>0$.
Williamson's theorem states that there exists a symplectic matrix $S$ such that
\[
        S^{\mathsf T}MS=
        \begin{pmatrix}
        \Lambda & 0\\
        0 & \Lambda
        \end{pmatrix},
        \qquad
        \Lambda=\operatorname{diag}(\lambda_1,\dots,\lambda_n),
\]
where the $\lambda_j$ are the symplectic eigenvalues \cite{Williamson1936,deGosson2006}.
For ellipsoids, all normalized symplectic capacities coincide and are given by
\[
        c(\Omega_M)=\frac{\pi}{\lambda_{\max}},
\]
where $\lambda_{\max}=\max_j\lambda_j$ \cite{HoferZehnder1994,deGosson2006}.

This formula is useful for linearized numerical tests because a ball evolved by a linear symplectic map becomes an ellipsoid.
It should not, however, be confused with a capacity formula for an arbitrary fixed-energy reaction surface.
A fixed-energy surface is odd-dimensional and unbounded in the reaction direction.
To obtain a finite capacity question one must specify a bounded full-dimensional domain, such as a ball, an ellipsoid, or a localized energy layer near the bottleneck.

\subsection{Connection to classical action}

In one degree of freedom, the symplectic form is the ordinary area form on the phase plane.
For connected simply connected planar domains, symplectic capacity reduces to area \cite{Siburg1993,deGosson2006}.
For an oscillator, the action
\[
        J=\frac{1}{2\pi}\oint p\,dq
\]
is the area enclosed by a periodic orbit divided by $2\pi$.
Thus a disk $p^2+q^2\le 2J_{\max}$ in a canonical oscillator plane has area
\[
        2\pi J_{\max}.
\]

This observation is the bridge to normal-form reaction dynamics.
Near an index-1 saddle, the bath modes are represented by canonical oscillator pairs $(q_k,p_k)$ with actions
\[
        J_k=\frac12(p_k^2+q_k^2)
\]
in the quadratic approximation.
The NHIM bounds the dividing surface, and its projection into each bath plane determines maximal admissible bath actions.
The quantities $2\pi J_k^{\max}(E)$ are therefore natural bath-plane area scales.
The central question of this paper is how far these action-area scales can be used as candidate symplectic widths for appropriately defined full-dimensional reactive neighborhoods.

\section{Phase Space Structures for Saddle--Center--\dots--Center Equilibria}

In the previous section, we established that evaluating a finite symplectic 
capacity requires a bounded phase-space domain defined by canonically 
conjugate planes. To apply this abstract topological constraint to a 
physical reaction channel, we must map it onto the concrete geometric 
structures that govern the transition. In this section, we explicitly 
construct these structures for a multidimensional Hamiltonian system. We 
take as our fundamental model the saddle--center--\dots--center equilibrium, 
which represents the local phase-space model associated with an index-1 saddle 
on a potential energy surface.

To mathematically realize the relevant phase-space architecture---and thereby 
enable a rigorous analysis of the reaction dynamics---we employ 
Poincar\'e--Birkhoff normal form theory. This constructive procedure does 
not simply decouple the system; rather, it provides a locally integrable 
approximation where the normal form Hamiltonian becomes a function of \(n\) 
independent integrals of motion (one associated with the reaction coordinate 
and \(n-1\) with the bounded bath modes). Within these adapted normal-form 
coordinates, we can explicitly define the standard local geometric hierarchy of the 
transition state: the \((2n-2)\)-dimensional dividing surface, the invariant 
\((2n-3)\)-dimensional Normally Hyperbolic Invariant Manifold (NHIM) that 
serves as its equator, and the local directional flux. Together, these 
explicitly computable objects provide the local geometric framework required 
to apply the symplectic capacity constraints developed in Section 3.

\subsection{Equilibrium point and its linearization}

In the phase space formulation of reaction dynamics, the key geometric structures arise from equilibrium points of Hamilton's equations that are of \emph{saddle--center--\dots--center} stability type.
Such an equilibrium corresponds, for a standard Hamiltonian of the form ``kinetic energy + potential energy'' on an \(n\)-dimensional configuration space, to an index-1 saddle of the potential energy surface \cite{WigginsEtAl2025}.
After translating the equilibrium to the origin of the \(2n\)-dimensional phase space, the quadratic part of the Hamiltonian is
\[
H_2(z) = \frac{1}{2} z^{\mathsf T} D^2 H(0) \, z,
\]
where \(z \in \mathbb{R}^{2n}\) and \(D^2 H(0)\) is the \(2n \times 2n\) symmetric Hessian matrix evaluated at the origin.
The corresponding linearized Hamilton's equations, \(\dot{z} = J D^2 H(0) z\) (where \(J\) is the standard \(2n \times 2n\) symplectic matrix), possess exactly one pair of real eigenvalues \(\pm \lambda\) and \(n-1\) pairs of purely imaginary eigenvalues \(\pm i\omega_k\) for \(k=2,\dots,n\), with \(\lambda, \omega_k > 0\).
This specific spectral configuration is the defining signature of a saddle--center--\dots--center equilibrium.

\subsection{Poincar\'e--Birkhoff normal form}

To understand the nonlinear dynamics near such an equilibrium, one employs the \emph{Poincar\'e--Birkhoff normal form} (or simply ``classical normal form'').
This is a constructive algorithm that, through a sequence of symplectic transformations, systematically simplifies the Taylor expansion of the Hamiltonian degree by degree in the local phase-space coordinates.
The result is a new set of coordinates, often called ``normal form coordinates,'' in which the Hamiltonian takes the form
\[
H_{\mathrm{NF}}(\bar{q},\bar{p}) = H_2(\bar{q},\bar{p}) + \text{higher-order terms},
\]
where the quadratic part \(H_2\) is unchanged, but the higher-order terms are specifically arranged to commute with \(H_2\) in the sense of Poisson brackets.
For a non-resonant index-1 saddle---meaning the purely imaginary bath frequencies \(\omega_k\) are rationally independent such that \(\sum_{k=2}^n m_k \omega_k \neq 0\) for any sequence of integers \(m_k\) not all zero---this procedure yields a locally integrable system: the normal form Hamiltonian becomes a function of \(n\) independent integrals.
A detailed exposition can be found in \cite{Waalkens2008, Uzer2002, Deprit1969}.

\subsection{Integrals of the normal form}

In the normal form coordinates, after a suitable linear symplectic change that brings the quadratic part to the standard form
\[
H_2 = \lambda \bar{p}_1 \bar{q}_1 + \sum_{k=2}^n \frac{\omega_k}{2}(\bar{p}_k^2 + \bar{q}_k^2),
\]
the higher-order terms can be expressed as polynomials in the following integrals:
\[
I = \bar{p}_1 \bar{q}_1, \qquad
J_k = \frac{1}{2}(\bar{p}_k^2 + \bar{q}_k^2),\quad k=2,\dots,n.
\]
These are constants of the motion for the truncated normal form Hamiltonian, and they are approximately conserved for the full system near the equilibrium.
The quantity \(I\) is associated with the ``reaction coordinate''; it is positive for reactive trajectories and negative for non-reactive ones.
The \(J_k\) are the actions of the bath modes (the ``centre'' degrees of freedom).
Rather than implying a strict decoupling of the modes, their conservation reflects the local integrability of the normal form approximation, where the effective frequencies of the dynamics depend nonlinearly on the excitation of all modes \cite{Waalkens2008}.

\subsection{Energy surfaces and the bottleneck}

For a fixed energy \(E\) above the saddle energy (taken as \(E_0=0\) for simplicity), the energy surface \(H_{\mathrm{NF}}=E\) has, in a local neighbourhood of the equilibrium, the geometric structure of a product space \(S^{2n-2}\times\mathbb{R}\).
This represents the multidimensional ``bottleneck'': the energy surface narrows to a minimal waist near the saddle.
To explicitly construct the dividing surface within this bottleneck, it is analytically convenient to rotate the reaction plane.
We introduce a linear transformation in the \((\bar{q}_1,\bar{p}_1)\) plane:
\[
Q_1 = \frac{1}{\sqrt{2}}(\bar{q}_1 - \bar{p}_1),\qquad P_1 = \frac{1}{\sqrt{2}}(\bar{q}_1 + \bar{p}_1).
\]
It is crucial that this specific rotation is a \emph{symplectic} transformation (preserving the area form \(dP_1 \wedge dQ_1 = d\bar{p}_1 \wedge d\bar{q}_1\)).
Because it is symplectic, the new coordinates \((Q_1, P_1)\) remain canonically conjugate, Hamilton's equations retain their standard form, and the rigorous topological constraints of symplectic capacity established in Section 3 remain fully valid.

In these rotated coordinates, dropping the bars on the bath modes for notational simplicity, the quadratic part of the Hamiltonian becomes
\[
H_2 = \frac{\lambda}{2}(P_1^2 - Q_1^2) + \sum_{k=2}^n \frac{\omega_k}{2}(p_k^2+q_k^2).
\]
The \((2n-2)\)-dimensional dividing surface (DS) is naturally defined by the geometric condition \(Q_1 = 0\).
Its intersection with the \((2n-1)\)-dimensional energy surface,
\[
\mathrm{DS}(E) = \left\{ Q_1=0,\;
\frac{\lambda}{2} P_1^2 + \sum_{k=2}^n \omega_k J_k = E \right\},
\]
forms a \((2n-2)\)-dimensional sphere, \(S^{2n-2}\).
The hemisphere with \(P_1>0\) is denoted \(\mathrm{DS}_{\mathrm{in}}(E)\) and corresponds to reactive trajectories entering the bottleneck region (i.e., moving from reactants toward products).
The hemisphere with \(P_1<0\) is \(\mathrm{DS}_{\mathrm{out}}(E)\), corresponding to exiting trajectories.
These hemispheres are strictly transverse to the Hamiltonian flow everywhere except on their common equatorial boundary \cite{Uzer2002, Waalkens2004}.

\subsection{Normally hyperbolic invariant manifold (NHIM)}

The equator of the dividing surface, given by \(Q_1 = P_1 = 0\), is a \((2n-3)\)-dimensional sphere:
\[
\mathrm{NHIM}(E) = \bigl\{ Q_1=0,\;
P_1=0,\; \sum_{k=2}^n \omega_k J_k = E \bigr\}.
\]
It is invariant because the equations of motion in the rotated coordinates give \(\dot{Q}_1 = \lambda P_1\) and \(\dot{P}_1 = \lambda Q_1\);
if both are zero initially, they remain zero for all time.
Moreover, the linearised dynamics normal to this manifold (in the directions of \(Q_1\) and \(P_1\)) is hyperbolic: one direction expands exponentially and the other contracts, with rates \(\pm \lambda\).
The dynamics tangent to the manifold is described by the centre degrees of freedom, which are oscillatory with frequencies \(\omega_k\).
This combination of stability properties is precisely the definition of a \emph{normally hyperbolic invariant manifold} \cite{Wiggins1994}.
In the reaction dynamics terminology, the NHIM thus plays the role of the “activated complex” in phase space \cite{Wiggins2001, Uzer2002}.

\subsection{Flux and action space volume}

The magnitude of the flux across the dividing surface is given by the symplectic area of \(\mathrm{DS}_{\mathrm{in}}(E)\):
\[
\phi(E) = \int_{\mathrm{DS}_{\mathrm{in}}(E)} d\sigma,
\]
where \(d\sigma\) is the restriction of the symplectic \((2n-2)\)-form \(\omega^{n-1}/(n-1)!\) to the hemisphere.

In the rotated normal form coordinates, this integral simplifies dramatically.
Because the coordinates \((P_1,J_2,\dots,J_n,\phi_2,\dots,\phi_n)\) are symplectic (with \(\phi_k\) the angles conjugate to \(J_k\)), the measure \(d\sigma\) becomes \(dJ_2\cdots dJ_n\, d\phi_2\cdots d\phi_n\).
Integrating over the angles gives a factor \((2\pi)^{n-1}\), and the condition \(P_1 = \sqrt{2(E - \sum \omega_k J_k)/\lambda}\) comes from the energy constraint.
Hence
\[
\phi(E) = (2\pi)^{n-1} \mathcal{V}(E),\qquad
\mathcal{V}(E) = \int_{\sum_{k=2}^n \omega_k J_k \le E} dJ_2\cdots dJ_n.
\]
Here, \(\mathcal{V}(E)\) is the volume of the solid region in action space defined by the energy condition.
Geometrically, the boundary of this allowed action-space region corresponds exactly to the NHIM (where the energy is entirely contained in the bath modes, meaning \(P_1 = Q_1 = 0\)).
Therefore, the total directional flux is precisely proportional to the volume enclosed by the NHIM in the transverse action space.

For the quadratic (linearized) Hamiltonian, this boundary forms a flat simplex, yielding an exact analytical volume.
For the full anharmonic Hamiltonian, the normal form provides a locally integrable approximation, and the identical formula holds with the quadratic energy boundary replaced by the nonlinear normal form polynomial \(K_{\mathrm{CNF}}(0,J_2,\dots,J_n)=E\).
This geometric realization is crucial for the arguments that follow.
It demonstrates not only how the normal form makes the multidimensional flux calculation tractable and physically interpretable \cite{Waalkens2004, Waalkens2005a, Waalkens2005b}, but it also establishes the fundamental phase-space architecture: while the flux measures the \emph{total volume} bounded by the NHIM, the symplectic capacity constraints explored in Section 5 will be dictated by the \emph{transverse widths} (the maximal action intercepts) of this exact same bounding manifold.

\subsection{Summary of phase space structures}

The key geometric objects and their representations in the rotated normal form coordinates are summarized in Table~1.

\begin{table}[ht]
\centering
\caption{Phase space structures in rotated normal form coordinates (energy \(E>0\)).}
\begin{tabular}{|l|l|}
\hline
Geometrical structure & Equation in rotated coordinates \\
\hline
Dividing surface & \(Q_1 = 0\) \\
Forward hemisphere (\(\mathrm{DS}_{\mathrm{in}}\)) & \(Q_1 = 0,\ P_1 > 0\) \\
Backward hemisphere (\(\mathrm{DS}_{\mathrm{out}}\)) & \(Q_1 = 0,\ P_1 < 0\) \\
NHIM & \(Q_1 = P_1 = 0\) \\
Local stable manifold (linearized) & \(P_1 = -Q_1\) \\
Local unstable manifold (linearized) & \(P_1 = Q_1\) \\
\hline
\end{tabular}
\end{table}

These structures, computed via the normal form transformation, provide a complete local description of the reaction dynamics and are the foundation for the symplectic capacity analysis that follows.
They do not, by themselves, identify a unique bounded full-dimensional reactive domain whose exact symplectic capacity is already known;
rather, they motivate the candidate width scales introduced in the next section.

\section{A Candidate Symplectic Width for the Reactive Region}

The phase-space structures constructed in Section 4 identify the transition-state geometry near an index-1 saddle: a dividing surface, its NHIM equator, stable and unstable manifolds, and a flux determined by the bath-action region enclosed by the NHIM.
We now ask a different question.
Can this same normal-form geometry be used to define a symplectic width scale for the local reactive bottleneck?

The question is motivated by Gromov's theorem, but it is not automatically answered by it.
Non-squeezing applies to full-dimensional symplectic domains.
The fixed-energy surface $\Sigma_E$ and the dividing surface inside it are not themselves full-dimensional symplectic domains.
Therefore the correct object is not the energy surface alone, but a localized full-dimensional neighborhood of the bottleneck.
The candidate width scales derived below should be understood in this sense.
They are action-area scales associated with a bounded, thickened proxy for the reactive channel.

\subsection{Why a full-dimensional proxy is needed}

The exact energy surface
\[
        \Sigma_E=\{z:H(z)=E\}
\]
is $(2n-1)$-dimensional and is not a symplectic manifold in its own right.
Moreover, the local energy surface near a saddle is unbounded in the reaction direction.
Thus assigning a finite symplectic capacity directly to $\Sigma_E$ is not well posed.
To formulate a capacity question one must specify a full-dimensional subset of the ambient phase space.

One natural construction is a narrow energy layer
\[
        \mathcal L_{E,\Delta E}
        =
        \{z:E-\Delta E\le H(z)\le E+\Delta E\}.
\]
Here $E$ is the central energy of the layer.
The layer is full-dimensional, but still too large unless it is localized near the transition-state region.
A dynamical localization can be made using gap-time or flow-time coordinates.
Near a dividing surface, points may be parametrized by coordinates on the dividing surface, an energy coordinate, and a time coordinate along the Hamiltonian flow.
This is the local product structure underlying gap-time parametrizations in phase-space reaction dynamics \cite{Ezra2009}.

Accordingly, a bounded local proxy for the bottleneck may be viewed schematically as
\[
        \mathcal N_{E,\Delta E,\Delta\psi}
        \approx
        \{\text{points launched from }\mathrm{DS}_{\mathrm{in}}(E'):
        E'\in[E-\Delta E,E+\Delta E],\ 0\le \psi\le\Delta\psi\},
\]
where $\psi$ is the flow-time coordinate.
The details of this thickening matter for an exact capacity calculation.
For this reason we do not claim that the formulas below are already a theorem giving the Gromov width of $\mathcal N_{E,\Delta E,\Delta\psi}$.
Rather, the layer construction explains why a full-dimensional domain is required and why the central energy $E$ should enter the leading-order bottleneck scale.

For a sufficiently narrow layer, the variation of the transverse bottleneck scale across the layer is small.
Thus the leading-order comparison is made at the central energy $E$.
If the narrowest point of the layer is needed, it occurs at $E-\Delta E$; in the limit $\Delta E\to0$, the central-energy and lower-edge scales coincide.
This is the sense in which the quantities $J_k^{\max}(E)$ and $c_{\mathrm{cand}}(E)$ below should be interpreted.

\subsection{Candidate width for the quadratic saddle}

For the quadratic saddle--center--$\cdots$--center Hamiltonian in rotated normal-form coordinates,
\[
        H_2
        =
        \frac{\lambda}{2}(P_1^2-Q_1^2)
        +
        \sum_{k=2}^n \omega_k J_k,
        \qquad
        J_k=\frac12(p_k^2+q_k^2),
\]
the forward hemisphere of the dividing surface is
\[
        \mathrm{DS}_{\mathrm{in}}(E)
        =
        \left\{
        Q_1=0,
        \ P_1>0,
        \ \frac{\lambda}{2}P_1^2+
        \sum_{k=2}^n\omega_kJ_k=E
        \right\}.
\]
The allowed bath actions satisfy
\[
        \sum_{k=2}^n\omega_kJ_k\le E,
        \qquad
        J_k\ge0.
\]
The NHIM is the boundary of this hemisphere obtained by setting $Q_1=P_1=0$; on it all the energy is stored in the bath actions.

For two degrees of freedom, there is one bath action and
\[
        0\le J_2\le \frac{E}{\omega_2}.
\]
The projection of the NHIM-bounded bath region onto the canonical plane $(q_2,p_2)$ is the disk
\[
        p_2^2+q_2^2\le 2\frac{E}{\omega_2}.
\]
The area of this disk is
\[
        2\pi\frac{E}{\omega_2}=2\pi J_2^{\max}(E).
\]
This is the exact transverse bath-plane area scale of the quadratic bottleneck.
It is also the capacity of the corresponding canonical cylinder based on that bath plane.

For three degrees of freedom, the actions satisfy
\[
        \omega_2J_2+
        \omega_3J_3\le E,
        \qquad
        J_2,J_3\ge0.
\]
The two canonical bath-plane area scales are
\[
        2\pi\frac{E}{\omega_2},
        \qquad
        2\pi\frac{E}{\omega_3}.
\]
The limiting transverse width is the smaller of these:
\[
        c_{\mathrm{cand}}(E)
        =
        2\pi\min\left(\frac{E}{\omega_2},\frac{E}{\omega_3}\right)
        =
        \frac{2\pi E}{\max(\omega_2,\omega_3)}.
\]
More generally, for the $n$-degree-of-freedom quadratic model,
\[
        J_k^{\max}(E)=\frac{E}{\omega_k},
        \qquad k=2,\dots,n,
\]
and the natural candidate width scale is
\[
        c_{\mathrm{cand}}(E)
        =
        2\pi\min_{k\ge2}J_k^{\max}(E)
        =
        2\pi\min_{k\ge2}\frac{E}{\omega_k}.
\]

In the quadratic model this formula has a precise geometric meaning: it is the smallest canonical bath-plane cylinder scale determined by the NHIM at energy $E$.
For the thickened local domain, it gives the leading-order central-energy width against which incoming full-dimensional sets should be compared.

\subsection{Anharmonic normal forms and a candidate reactive width}

For the Eckart--Morse and Eckart--Morse--Morse models, the full dynamics is anharmonic and not globally integrable.
Nevertheless, the high-order classical normal form gives a locally integrable approximation near the saddle.
In this approximation the Hamiltonian can be written as a polynomial
\[
        K_{\mathrm{CNF}}(I,J_2,\dots,J_n)
        =
        E_0+\lambda I+
        \sum_{k=2}^n\omega_kJ_k+
        \text{higher-order terms},
\]
where $I$ is the reaction integral and the $J_k$ are bath actions.
On the NHIM, $Q_1=P_1=0$, so $I=0$.
Thus the bath-action boundary of the bottleneck at energy $E$ is determined by
\[
        K_{\mathrm{CNF}}(0,J_2,\dots,J_n)=E,
        \qquad
        J_k\ge0.
\]

For each bath mode we define
\[
        J_k^{\max}(E)
        =
        \sup\left\{
        J_k\ge0:
        \exists (J_2,\dots,J_n)\text{ satisfying }
        K_{\mathrm{CNF}}(0,J_2,\dots,J_n)=E
        \right\}.
\]
Equivalently, in practical computations one often obtains $J_k^{\max}(E)$ by setting the other bath actions to zero and solving the resulting one-dimensional polynomial equation, provided the local normal-form energy boundary is well behaved.
The candidate anharmonic width scale is then
\[
        c_{\mathrm{cand}}(E)
        =
        2\pi\min_{k\ge2}J_k^{\max}(E).
\]

For example, in two degrees of freedom suppose
\[
        K_{\mathrm{CNF}}(I,J_2)
        \approx
        E_0+
        \lambda I+
        \omega_2J_2+
        \alpha J_2^2+
        \beta I J_2.
\]
On the NHIM, $I=0$, and $J_2^{\max}$ is determined by
\[
        \omega_2J_2+
        \alpha J_2^2=E-E_0.
\]
The positive root is
\[
        J_2^{\max}(E)
        =
        \frac{-\omega_2+
        \sqrt{\omega_2^2+4\alpha(E-E_0)}}{2\alpha},
\]
which reduces to $(E-E_0)/\omega_2$ in the limit $\alpha\to0$.
The corresponding candidate width is $2\pi J_2^{\max}(E)$.

\begin{conjecture}[Candidate reactive width for anharmonic normal forms]
For energies sufficiently close to an index-1 saddle, the normal-form geometry of the NHIM suggests the transverse reactive width scale
\[
        c_{\mathrm{cand}}(E)
        =
        2\pi\min_{k\ge2}J_k^{\max}(E),
\]
where $J_k^{\max}(E)$ is obtained from the truncated classical normal form evaluated on the NHIM, $I=0$.
This quantity should be interpreted as a candidate width for a bounded full-dimensional neighborhood of the bottleneck, not as an already proven formula for the exact Gromov width of an arbitrary reactive domain.
\end{conjecture}

This cautious formulation is important.
The normal form provides canonical coordinates and computable bath-action boundaries.
It does not by itself solve the global embedding problem for a chosen thickened domain.
The candidate width scale is therefore a diagnostic and organizing quantity: it identifies the smallest transverse canonical bath-plane area scale suggested by the NHIM geometry.

\subsection{Why introduce a symplectic width beyond the flux?}

The directional flux measures the total phase-space volume crossing a dividing surface per unit time.
It is a central quantity in transition state theory and remains indispensable.
A symplectic-width viewpoint asks a different question: how is the transported set arranged in canonical transverse directions?
Two full-dimensional ensembles may have comparable volume, comparable energy spread, or comparable nominal flux, but very different two-dimensional projections onto the bath-action planes.

This distinction is invisible to Liouville volume alone.
It is precisely the distinction emphasized by non-squeezing: Hamiltonian flow cannot treat canonical coordinate--momentum directions as arbitrary compressible directions.
The quantities $2\pi J_k^{\max}(E)$ therefore provide a natural set of transverse scales against which incoming ensembles may be compared.
If an ensemble is strongly localized near a high-action boundary in one bath mode, it may have little reactive action available and may display substantial finite-time delay, even though the total energy lies above the saddle.

The numerical experiments below should be read in this diagnostic sense.
They do not prove that $c_{\mathrm{cand}}(E)$ is a universal sharp threshold for reaction.
They illustrate how bath-action geometry, finite-time transmission, and symplectic width scales can interact in concrete normal-form models.

\section{High-Order Normal Forms: Eckart--Morse Models}\label{sec:highorder}

In Section 5, we introduced a geometrically motivated candidate width scale 
based on the maximal transverse bath actions. However, the purely quadratic 
models used to derive those initial area scales lack the fundamental features 
of realistic molecular reactions: anharmonicity and nonlinear mode coupling. 
To determine if the proposed geometric constraints actually govern transport 
in physically realistic phase spaces, we must transition from exactly solvable 
linear models to fully anharmonic Hamiltonians.

In this section, we introduce the Eckart--Morse (2 DoF) and Eckart--Morse--Morse 
(3 DoF) models. These systems serve as standard paradigms in multidimensional 
reaction dynamics, accurately modeling an asymmetric reaction barrier coupled 
to anharmonic transverse oscillations. Because these full systems are 
non-integrable, we cannot compute their bounding phase-space structures 
analytically. Instead, we utilize high-order Poincar\'e--Birkhoff normal 
forms to construct a locally integrable approximation near the index-1 saddle. 
This mathematical machinery allows us to explicitly extract the relevant 
geometric quantities---the total directional flux and the maximal admissible 
bath actions---providing the exact computable parameters required to test the 
symplectic-width hypothesis in the numerical experiments of Section 7.

\subsection{The Eckart--Morse Hamiltonian}

To specify the model used in the numerical experiments and to ground the 
normal-form analysis in a concrete physical model, we explicitly define the 
Hamiltonian for the Eckart--Morse (2 DoF) and Eckart--Morse--Morse (3 DoF) 
systems. These models serve as standard paradigms in reaction dynamics, 
accurately representing an asymmetric reaction channel coupled to anharmonic 
transverse vibrational modes \cite{Waalkens2008}.

For the full 3 DoF Eckart--Morse--Morse system, the Hamiltonian in physical 
phase-space coordinates \((x, y, z, p_x, p_y, p_z)\) is given by
$$H = \frac{1}{2m} (p_x^2 + p_y^2 + p_z^2) + V_E(x) + V_{M,y}(y) + V_{M,z}(z) + \epsilon (p_x p_y + p_x p_z + p_y p_z).$$
Here, the reaction coordinate \(x\) is governed by the Eckart potential:
$$V_E(x) = A \frac{e^{(x+x_0)/a}}{1+e^{(x+x_0)/a}} + B \frac{e^{(x+x_0)/a}}{(1+e^{(x+x_0)/a})^2},$$
where \(A\) controls the asymptotic asymmetry (the endo/exothermicity of the 
reaction), \(B\) scales the barrier height, \(a\) is the width parameter, and 
\(x_0\) is a spatial shift.

The transverse bath modes \(y\) and \(z\) are modeled by Morse oscillators, 
representing the chemical bonds:
$$V_{M,k}(q) = D_e \left( e^{-2a_M q} - 2e^{-a_M q} \right), \qquad (q = y, z),$$
where \(D_e\) is the dissociation depth and \(a_M\) controls the anharmonicity 
of the well.

The cross-derivative terms parameterized by \(\epsilon\) in the kinetic energy 
induce momentum-coupling between the reaction coordinate and the bath modes, 
as well as between the bath modes themselves. (The 2 DoF Eckart--Morse system 
is simply the reduction of this Hamiltonian to the \((x,y,p_x,p_y)\) subspace, 
dropping the \(z\)-dependent terms).

The specific numerical parameters used in our finite-time transmission 
experiments are taken directly from the benchmark study by Waalkens 
\emph{et al.}~\cite{Waalkens2008} (e.g., \(m=1\), \(\epsilon=0.3\)). Because 
this full Hamiltonian is non-integrable due to the kinetic coupling and 
the anharmonic potentials, one cannot compute the exact bounding phase-space 
structures analytically. Instead, we compute its high-order Poincar\'e--Birkhoff 
normal form to construct the locally integrable approximation required to 
extract the symplectic invariants.

\subsection{Extracting geometric quantities from the normal form}

Once the high-order normal form polynomial $K_{\mathrm{CNF}}(I, J_2, \dots, J_n)$ 
has been computed for these systems, we can extract its geometric invariants. 
Because the normal form effectively straightens out the local dynamics, it 
allows us to evaluate action-space volumes and candidate width scales for 
the thickened proxy neighborhoods near the bottleneck semi-analytically.

Consistent with our approach in Section 5.1, we evaluate the transverse 
geometry of the thickened domain at its central target energy $E$. We isolate 
the bottleneck cross-section by restricting the dynamics to the invariant 
subspace corresponding to $I=0$ (the dividing surface). By analyzing the 
truncated polynomial $K_{\mathrm{CNF}}(0, J_2, \dots, J_n) = E$, we compute 
three primary quantities:

\begin{itemize}
    \item \textbf{The directional flux} $\phi(E)$: The flux is proportional 
    to the Euclidean volume of the corresponding positive-action region, given 
    by $\phi(E) = (2\pi)^{n-1} V(E)$. Here, $V(E)$ is the volume in the 
    positive action space bounded by the contour at the central energy:
    $$ \Omega_E = \left\{ (J_2, \dots, J_n) \in \mathbb{R}_{\geq 0}^{n-1} \mid K_{\mathrm{CNF}}(0, J_2, \dots, J_n) \leq E \right\}. $$
    Because $K_{\mathrm{CNF}}$ is a non-linear, high-degree polynomial, an 
    analytical volume formula is generally unavailable. Instead, $V(E)$ is 
    evaluated via numerical integration. Depending on the dimensionality $n$, 
    we employ Monte Carlo integration, utilizing the maximal actions (defined 
    below) to construct an efficient bounding box for the sampling domain.

    \item \textbf{The maximal actions} $J_k^{\max}(E)$: Geometrically, these 
    represent the maximum allowable limits of the physically accessible region 
    in the transverse action space for the central target energy $E$. For each 
    degree of freedom $k \in \{2, \dots, n\}$, we determine $J_k^{\max}(E)$ by 
    setting all other transverse actions to zero and finding the positive real 
    root of the univariate polynomial equation:
    $$ K_{\mathrm{CNF}}(0, \dots, 0, J_k, 0, \dots, 0) = E. $$
    Assuming the normal form retains the definiteness of the underlying 
    Hamiltonian near the equilibrium, this equation yields a unique positive 
    root. This root can be rapidly located using standard numerical root-finding 
    algorithms (e.g., Newton-Raphson or Brent's method).

    \item \textbf{The candidate width scale} $c_{\mathrm{cand}}(E)$: In the 
    purely quadratic case, the candidate width of the thickened proxy domain 
    is governed entirely by the minimal maximal action. By direct analogy, 
    we extend this candidate width scale to the non-linear regime captured 
    by our high-order model, defined as:
    $$ c_{\mathrm{cand}}(E) = 2\pi \min_{k} J_k^{\max}(E). $$
\end{itemize}

By calculating these three quantities across a range of central target 
energies $E$, we can examine how the candidate width scale varies relative 
to the flux and the maximal admissible bath actions in the nonlinear normal form.

\section{Numerical Experiments}

We now present numerical experiments designed to probe the geometric ideas developed above.
The goal is not to turn Gromov's theorem into a direct formula for a reaction probability.
Rather, the goal is to separate three related diagnostics:
\begin{enumerate}[label=(\roman*)]
    \item the invariant symplectic size of a full-dimensional set;
    \item the ordinary two-dimensional projection of that set onto selected canonical planes;
    \item the finite-time transmission behavior of trajectories drawn from the set.
\end{enumerate}
The first is controlled by symplectic topology, the second is a visible but non-invariant measure of deformation, and the third depends on the actual Hamiltonian dynamics and the chosen observation time.

The experiments use the normal-form framework developed in Sections 5 and 6.
In each case the relevant energy $E$ is the central energy of the ensemble or layer being tested.
The candidate width $c_{\mathrm{cand}}(E)$ is therefore a central-energy diagnostic of the bottleneck geometry, not a claim that a fixed-energy surface by itself has a symplectic capacity.

\subsection{Experiment 1: Linearized non-squeezing consistency check}

The first experiment is a controlled linear test.
A small four-dimensional ball is mapped by a symplectic linear transformation and then evolved under the linearized Hamiltonian flow.
Because the evolution is linear and symplectic, the image of the ball is an ellipsoid whose capacity can be computed exactly from its symplectic spectrum.
This setting is useful because it avoids filamentation and isolates the distinction between capacity and projected area.

To make the projection dynamics nontrivial, we first apply a symplectic mixing matrix $S_{\mathrm{mix}}$ to the initial ball.
The resulting ellipsoid is skewed relative to the normal-form coordinates.
We then evolve it backward in time using the exact state transition matrix $\Phi(-\tau)$ of the linearized Hamiltonian.
For comparison with the saddle dynamics, we monitor the ordinary projection area of the evolving ellipsoid onto the saddle plane $(Q_1,P_1)$:
\[
 A(\tau)
 =
 \pi r^2
 \sqrt{\det\left(
 P\,\Phi(-\tau)S_{\mathrm{mix}}S_{\mathrm{mix}}^{\mathsf T}
 \Phi(-\tau)^{\mathsf T}P^{\mathsf T}
 \right)},
\]
where $P$ is the $2\times4$ projection matrix onto the saddle coordinates.
This is an ordinary projection-area diagnostic.
In this special linear ellipsoidal setting, its minimum can be compared directly with the exact capacity scale $\pi r^2$.

\begin{figure}[htbp]
    \centering
    \includegraphics[width=\linewidth]{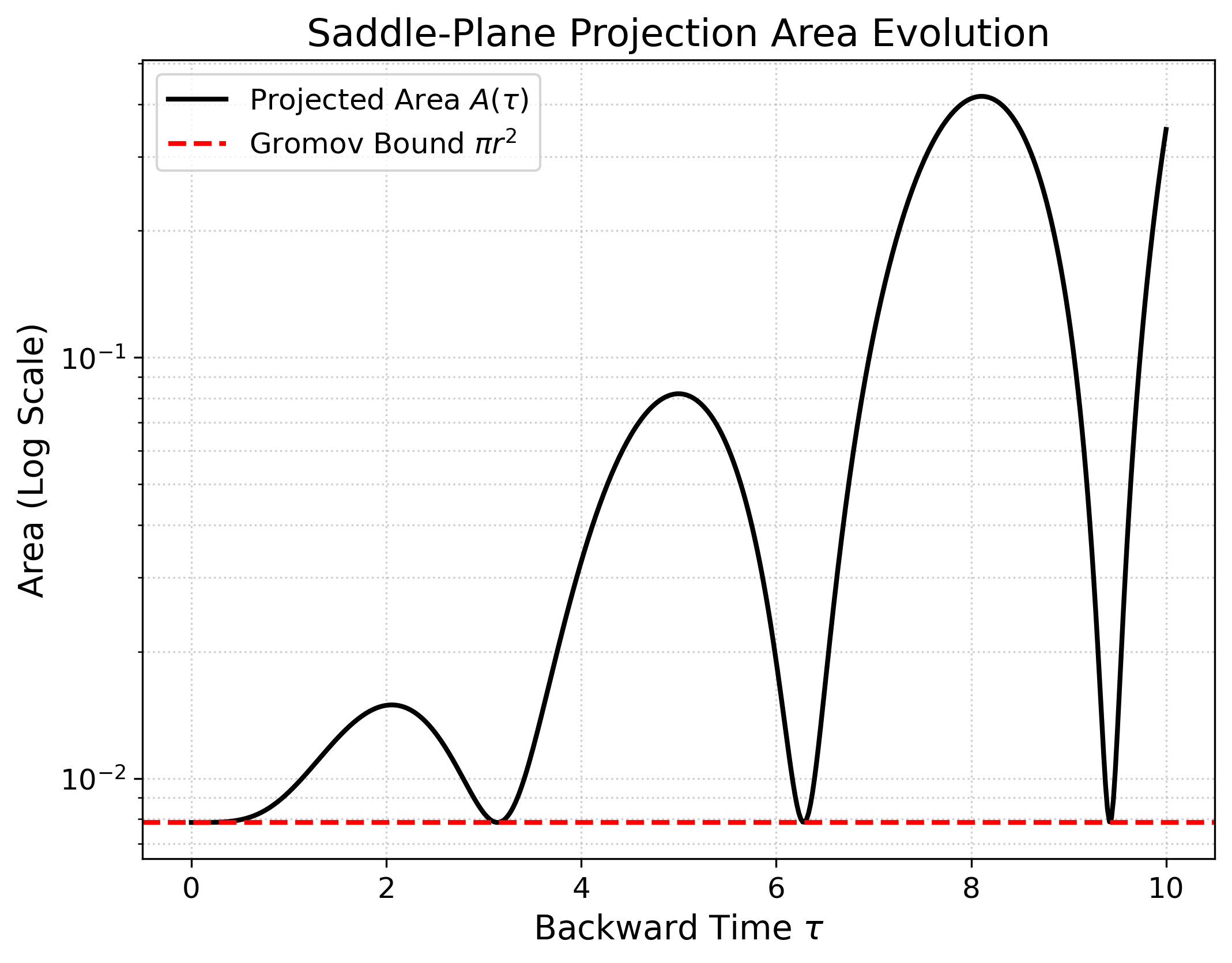}
    \caption{Evolution of the saddle-plane projection area $A(\tau)$ for a locally mixed four-dimensional phase-space ball evolved backward in time under the linearized Hamiltonian flow.  The dashed line marks the ball capacity $\pi r^2$.  In this linear ellipsoidal test, the projected area fluctuates as the ellipsoid tumbles but does not fall below the Gromov scale.}
    \label{fig:projection_area}
\end{figure}

Figure~\ref{fig:projection_area} shows the projected area during the backward evolution.
The area grows because of the hyperbolic stretching and oscillates because the skewed ellipsoid tumbles relative to the saddle plane.
The key observation is that the projection-area diagnostic touches but does not cross the Gromov scale $\pi r^2$.
This is a consistency check of the linear symplectic geometry, not a numerical proof of the theorem itself.

\begin{figure}[htbp]
    \centering
    \includegraphics[width=\linewidth]{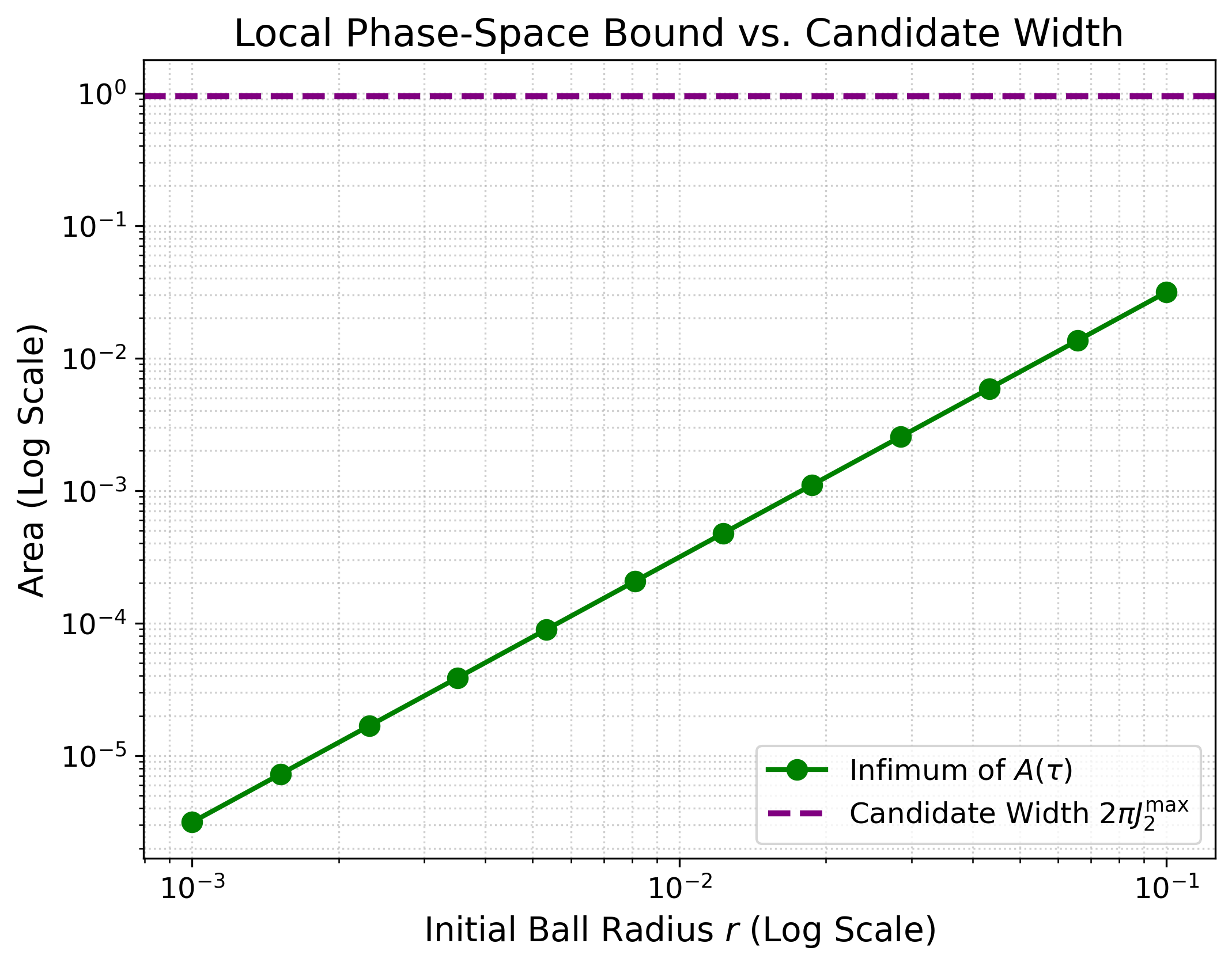}
    \caption{Minimum saddle-plane projection area $\inf_\tau A(\tau)$ for a range of initial ball radii $r$.  The observed scaling follows the local ball capacity $\pi r^2$.  The dashed reference level indicates the larger candidate bottleneck width $c_{\mathrm{cand}}(E)=2\pi J_2^{\max}(E)$ associated with the central energy of the normal-form bottleneck.}
    \label{fig:infimum_scaling}
\end{figure}

Figure~\ref{fig:infimum_scaling} repeats this calculation for several initial radii and plots the minimum projected area over the time interval considered.
The scaling with $\pi r^2$ separates the local capacity scale of the test ball from the larger candidate bottleneck width $c_{\mathrm{cand}}(E)$.

\subsection{Experiment 2: A bath-localized ensemble calculation}

The second experiment studies finite-time transmission for two ensembles in a two-degree-of-freedom anharmonic normal-form model.
The purpose is to test whether localization near a high bath-action boundary can produce a strong dynamical delay.
This is a trajectory diagnostic motivated by the symplectic-width picture; it is not a direct theorem-level consequence of non-squeezing.

This comparison should not be read as proving that the observed delay is invisible to all possible flux- or volume-based descriptions.
A sharper test of that claim would require constructing ensembles matched in total phase-space volume, flux, or microcanonical measure while differing primarily in their transverse bath-action localization.
The present calculation has a more limited purpose: it tests whether bath-action localization can produce a substantial finite-time transmission delay in a normal-form model, as suggested by the symplectic-width viewpoint.

Both ensembles contain $N=5000$ initial conditions in the forward-reactive half-space, with $Q_1<0$ and $P_1>0$, and both span a narrow energy interval $[E-\Delta E,E+\Delta E]$ about the same central energy $E$.
They differ in how the available action is distributed:
\begin{itemize}
    \item \textbf{Ensemble A (unbiased action-space sampling).}
    Initial conditions are sampled broadly across the physically accessible action region near the bottleneck.
    This ensemble provides a baseline for finite-time transmission.
    The term ``unbiased'' is meant geometrically; it does not assume thermal equilibrium.

    \item \textbf{Ensemble B (bath-localized sampling).}
    The bath action is restricted to a high-action wedge
    \[
        J_2\in[\xi J_2^{\max},J_2^{\max}],
    \]
    where $\xi\in[0,1]$.
    Increasing $\xi$ places more of the available energy budget into the transverse bath motion and leaves less room for reactive action.
\end{itemize}

For the normal-form truncation used here, the effective hyperbolic rate depends on bath action,
\[
        \Lambda(J_2)=\lambda+b_2J_2,
\]
with $\lambda=0.7350$ and $b_2=-0.0123$.
The reaction coordinate evolves according to
\[
        Q_1(t)=Q_1\cosh(\Lambda t)+P_1\sinh(\Lambda t).
\]
A trajectory is counted as transmitted if $Q_1(t)>0$ at the final observation time $t_{\mathrm{max}}=5/\lambda$.
Thus the plotted quantity is a finite-time transmission fraction, not an asymptotic reaction probability.

The bath-localized ensemble displays a strong reduction in finite-time transmission as $\xi$ increases, reaching zero transmission within the selected observation window for $\xi\gtrsim0.6$ in the reported calculation.
This should be interpreted as a finite-time delay effect.
Because $\Lambda(J_2)$ remains positive over the action range considered, the trajectories are not permanently trapped by this normal-form model.
They are slowed because high bath action changes the effective saddle rate and places the initial conditions closer, in practical finite-time terms, to the transition-state boundary between rapid crossing and delayed crossing.

\begin{figure}[htbp]
    \centering
    \includegraphics[width=0.8\textwidth]{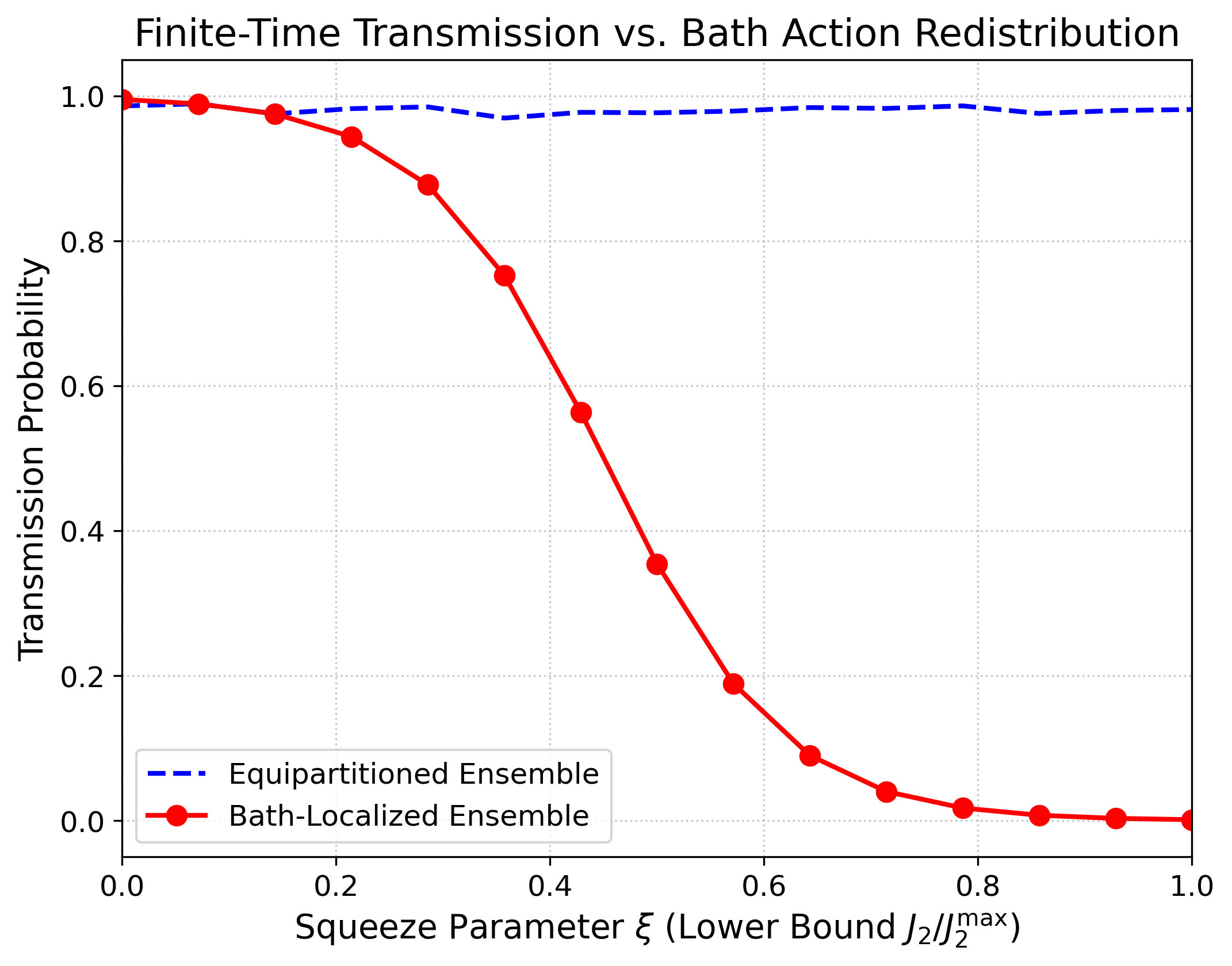}
    \caption{Finite-time transmission as a function of the bath-localization parameter $\xi$.  Both ensembles contain $N=5000$ trajectories prepared near the same central energy $E$.  The broadly sampled ensemble remains near full transmission in the chosen observation window, while the bath-localized ensemble displays a strong finite-time delay as its support is pushed toward $J_2^{\max}$.  The figure should be read as a trajectory diagnostic motivated by the candidate symplectic-width scale, not as a direct proof of a universal reaction threshold.}
    \label{fig:blockade_scan}
\end{figure}

This experiment isolates one possible dynamical manifestation of the symplectic-filter viewpoint.
The candidate width scale identifies the high bath-action boundary of the bottleneck.
When the initial ensemble is biased toward that boundary, finite-time transmission can be strongly suppressed.
The result is consistent with the idea that transverse bath-action geometry matters for transport, but the conclusion remains model- and time-window-dependent. 

\subsection{Implementation notes}

The required normal-form software is described in \cite{Waalkens2008} and can be downloaded at
\url{https://github.com/stephenwiggins/NormalForm}.
The linearized test uses the state transition matrix obtained from the quadratic normal form.
The bath-localized calculation uses the leading anharmonic coupling extracted from the classical limit of the tenth-order Eckart--Morse normal form.
The figures should be regenerated with a symplectic or exactly symplectic linear propagator where applicable, and for nonlinear trajectory integrations one should monitor both energy error and symplecticity error of the computed time-$t$ map.

\section{Discussion and Conclusion}

This paper has developed a symplectic-geometric perspective on classical reaction dynamics near an index-1 saddle.
The central proposal is that a reaction bottleneck can be viewed as a symplectic filter: a local transition-state structure whose transverse canonical width constrains how full-dimensional phase-space sets can pass through it.
The normal-form machinery supplies the reaction coordinates, the bath coordinates, the dividing surface, the NHIM, and the bath-action geometry needed to make this statement concrete.
Gromov's non-squeezing theorem supplies the underlying topological principle: a Hamiltonian image of a ball cannot be coherently squeezed into a smaller canonical cylinder.

The most important clarification is that fixed-energy transition-state geometry and symplectic capacity live in different categories.
A fixed-energy surface is odd-dimensional and is not the direct object of a capacity calculation.
To pose a non-squeezing question one must use a full-dimensional domain, such as a ball, an ellipsoid, or a narrow energy layer localized near the bottleneck.
The normal-form formulas in this paper should therefore be interpreted as central-energy candidate width scales for such localized full-dimensional domains.
They are not asserted to be exact capacities of the unthickened energy surface.

For quadratic saddle--center--$\cdots$--center models, the NHIM geometry gives explicit bath-action area scales.
For each transverse mode, the maximal action $J_k^{\max}(E)$ determines a canonical bath-plane disk of area $2\pi J_k^{\max}(E)$.
The smallest of these areas gives the natural candidate width
\[
        c_{\mathrm{cand}}(E)
        =
        2\pi\min_{k\ge2}J_k^{\max}(E).
\]
For anharmonic normal forms, the same construction can be carried out by solving the truncated normal-form equation on the NHIM, $I=0$.
This yields a computable diagnostic of how anharmonicity and mode coupling distort the transverse bottleneck geometry.

The numerical experiments should be read in this same diagnostic spirit.
The linear experiment checks that the projected-area behavior of a symplectically evolved ellipsoid is consistent with the local capacity scale of the initial ball.
The bath-localized experiment shows that placing an ensemble near a high bath-action boundary can produce strong finite-time delay in an anharmonic normal-form model.
Neither experiment proves a universal reaction-rate theorem.
Rather, they illustrate how the symplectic-width viewpoint can guide the choice of trajectory diagnostics: one should monitor not only flux and energy, but also the evolving projections of finite ensembles onto canonical bath planes, the distribution of bath actions, finite-time transmission, and possible returns.

Several mathematical problems remain open.
The first is to define, for a given reaction problem, a canonical full-dimensional reactive neighborhood whose symplectic capacity can be computed or estimated sharply.
The second is to determine when the normal-form action scale $c_{\mathrm{cand}}(E)$ coincides with a genuine Gromov width or cylindrical capacity of that neighborhood.
The third is to understand how these local capacity constraints interact with global nonlinear dynamics, including action exchange, chaotic scattering, return crossings, and the breakdown of the truncated normal-form approximation away from the saddle.

The main contribution of the paper is therefore not the proof of a new exact capacity theorem for chemical reactions.
It is the formulation of a concrete geometric framework linking NHIM-based transition-state geometry, bath-action width scales, full-dimensional energy layers, and finite-time trajectory diagnostics.
This framework suggests that mode-specific excitation should not be viewed only through the lens of total energy or total flux.
It also has a symplectic geometry: how the ensemble occupies canonical transverse directions can matter for how it negotiates the bottleneck.

\section{Computing Symplectic Eigenvalues}

This algebraic procedure is used to compute the symplectic capacity of the 
phase-space ellipsoids arising in the linearized test of Section 7.1, where 
the ellipsoid is represented by a positive-definite shape matrix. Let \(M\) be a positive 
definite symmetric \(2n\times 2n\) matrix. The symplectic eigenvalues of \(M\) 
are defined as the positive numbers \(\lambda_1 \ge \lambda_2 \ge \dots \ge \lambda_n > 0\) 
such that the eigenvalues of \(J M\) are \(\pm i \lambda_j\) (see \cite{Williamson1936}). 
They can be computed as follows.

\begin{enumerate}
    \item Compute the matrix \(W = M^{1/2} J M^{1/2}\). Because \(M\) is 
    symmetric positive definite, its square root \(M^{1/2}\) exists and is 
    symmetric. The matrix \(W\) is skew-symmetric, since \(J^{\mathsf T} = -J\) 
    and \(M^{1/2}\) is symmetric:
    $$ W^{\mathsf T} = (M^{1/2})^{\mathsf T} J^{\mathsf T} (M^{1/2})^{\mathsf T} = M^{1/2} (-J) M^{1/2} = -W. $$
    Hence \(W\) is real skew-symmetric, and its eigenvalues come in purely 
    imaginary conjugate pairs \(\pm i \lambda_j\) with \(\lambda_j \ge 0\).

    \item The numbers \(\lambda_j\) are precisely the symplectic eigenvalues 
    of \(M\) \cite{deGosson2006}. In fact, the eigenvalues of \(J M\) are the 
    same as those of \(M^{1/2} J M^{1/2}\) because \(M^{1/2} J M^{1/2} = M^{1/2} (J M) M^{-1/2}\) 
    is similar to \(J M\).

    \item Thus, one can compute the symplectic eigenvalues by diagonalising 
    the real skew-symmetric matrix \(W\) (e.g., using the Schur decomposition 
    or a standard eigenvalue routine) and taking the absolute values of the 
    imaginary parts of its eigenvalues. These are the \(\lambda_j\).
\end{enumerate}

For a block-diagonal matrix 
$$ M = \begin{pmatrix} A & 0 \\ 0 & B \end{pmatrix} $$ 
with \(A,B\) symmetric positive definite, one may also compute the symplectic 
eigenvalues as the square roots of the eigenvalues of \(AB\) \cite{HoferZehnder1994}. 
This follows from the fact that
$$ J M = \begin{pmatrix} 0 & I \\ -I & 0 \end{pmatrix} \begin{pmatrix} A & 0 \\ 0 & B \end{pmatrix} = \begin{pmatrix} 0 & B \\ -A & 0 \end{pmatrix}, $$
and the eigenvalues of 
$$ \begin{pmatrix} 0 & B \\ -A & 0 \end{pmatrix} $$ 
are \(\pm i \sqrt{\mu}\) where \(\mu\) runs over the eigenvalues of \(AB\).

\section{Verification of Classical Normal Form Coefficients for Eckart--Morse Models}

The classical normal form coefficients used in Sections 5, 6, and 7 are 
extracted from the formal classical limit of the 10th-order quantum normal 
form computed by Waalkens \emph{et al.}~\cite{Waalkens2008} (discarding all 
terms from the quantum Weyl symbols containing a factor of \(\hbar\)). The 
complete set of high-order polynomial coefficients for both the 2 DoF and 
3 DoF models can be found explicitly tabulated in the Appendices of that 
reference. For readers wishing to generate these coefficients independently 
or adapt them for modified potential parameters, the coefficients are 
computed using algorithmic Lie transform perturbation theory. The recursive 
algorithms required for this procedure, which are readily implemented in 
standard symbolic algebra systems (e.g., Mathematica or SymPy), are detailed 
comprehensively within the same benchmark study.

For the 2 DoF Eckart--Morse system evaluated in Experiment 2, the relevant 
truncated polynomial required to determine the bottleneck cross-section at a 
central target energy \(E\) and the leading-order transverse dynamics is:
$$ K_{\mathrm{CNF}}(I,J_2) \approx E_0 + \lambda I + \omega_2 J_2 + b_2 I J_2. $$
The fundamental parameters defining the saddle and the linear frequencies are 
the saddle energy \(E_0 = -0.9875\), the linear Lyapunov exponent \(\lambda = 0.7350\), 
and the harmonic bath frequency \(\omega_2 = 1.8225\). The leading anharmonic 
cross-coupling between the reaction coordinate and the transverse mode is 
\(b_2 = -0.0123\).

This explicit truncation illustrates how the linear approximation of the 
maximal bath action \(J_2^{\max}(E) \approx (E - E_0)/\omega_2\) is obtained 
on the dividing surface (\(I=0\)), while the effective Lyapunov exponent 
\(\Lambda(J_2) = \partial K_{\mathrm{CNF}} / \partial I = \lambda + b_2 J_2\) 
simultaneously emerges to govern the dynamical delay. These values, along 
with the corresponding parameters for the full 3 DoF Eckart--Morse--Morse 
system (e.g., the second bath frequency \(\omega_3 = 1.267\)), provide the 
necessary algebraic foundation for the numerical experiments discussed in 
the main text.

\end{document}